\documentclass[12pt]{article}
\usepackage[T2A]{fontenc}
\usepackage[cp1251]{inputenc}
\usepackage[english,russian]{babel}

\usepackage[tbtags]{amsmath}
\usepackage{amsfonts,amssymb}

\def\ba#1{\begin{array}{#1}}
\def\ea{\end{array}}
\def\beq#1{\begin{equation}\label{#1}}
\def\eeq{\end{equation}}

\usepackage{textcomp}
\def\No{\textnumero}

\def\dst{\displaystyle}

\newcommand{\Set}[1]{\Omega_{#1}}
\newcommand{\MSet}[1]{\Omega_{(#1)}}
\newcommand{\Size}[2]{\nu_{#1}(#2)}
\newcommand{\MSize}[2]{\nu_{(#1)}(#2)}
\newcommand{\eSize}[1]{n_{#1}}
\newcommand{\eMSize}[1]{n_{(#1)}}
\newcommand{\Zp}{\mathbb{Z}_+}
\newcommand{\ZZ}{\mathbb{Z}}
\newcommand{\FF}{\mathbb{F}}
\newcommand{\eps}{\varepsilon}
\newcommand{\Sym}{Sym}
\newcommand{\const}{\mathrm{const}}

\newcommand{\Gminuse}{G^\circ} 
\newtheorem{theorem}{Теорема}[section]

\newtheorem{proposition}[theorem]{Предложение}
\newtheorem{Corollary}[theorem]{Следствие}

\newcommand{\proof}{\noindent {\it Доказательство}. }
\newcommand{\eop}{\hfill$\Box$}
\newcommand{\Fix}[1]{\mathrm{Fix}_{#1}} 
\newcommand{\MFix}[1]{\mathrm{Fix}_{(#1)}} %
\newcommand{\Stab}{\mathrm{Stab}} 
\newcommand{\avStab}[2]{\left\langle\Stab^{#1}\right\rangle_{#2}} 
\newcommand{\avStabSet}[1]{\avStab{}{#1}}
\newcommand{\avStabMSet}[1]{\avStab{}{(#1)}}

\title{Об асимптотически свободном действии групп подстановок на подмножествах и мультимножествах
\footnote{Настоящий текст представляет собой переработанный вариант статьи автора <<Об асимптотически свободном действии групп перестановок на подмножествах и мультимножествах>> (<<Дискретная математика>>, 2014, т.~26, вып.~3, с.~101--120). По стечению обстоятельств в печать попал
дорецензионный вариант.}}

\author{С.Ю.~Садов}
\date{}

\begin{document}
\maketitle

\begin{abstract} 
Пусть $G$ --- группа подстановок, действующая на конечном множестве $\Omega$ мощности $n$.  Для числа орбит индуцированного действия $G$ на множестве $\Omega_m$ всех подмножеств $\Omega$ мощности $m$ имеет место тривиальная двусторонняя оценка $|\Omega_m|/|G|\leq |\Omega_m/G|\leq |\Omega_m|$. В статье даны улучшения правого неравенства в терминах минимальной степени группы $G$ или минимальной степени ее подмножества, дополнение которого мало. Рассмотрены приложения к задачам перечисления точечных конфигураций в аффинном пространстве над конечным полем, графов и гиперграфов. С использованием известных результатов теории групп подстановок показано, что если $G$ --- произвольная 2-транзитивная группа, за исключением симметрической и знакопеременной, $m$ и $n$ велики и отношение $m/n$ отделено от $0$ и $1$, то $|\Omega_m/G|\approx |\Omega_m|/|G|$.

Аналогичные результаты верны для индуцированного действия $G$ на множестве $\Omega_{(m)}$ всех мультимножеств $\Omega$ веса $m$, если отношение $m/n$ не слишком близко к $0$.
\end{abstract}

\medskip
{\small
{\bf Ключевые слова: } группа подстановок, регулярные орбиты, асимптотически свободное действие, минимальная степень группы, средний размер стабилизатора,  перечисление графов, перечисление аффинных конфигураций, асимптотика числа орбит.
}



\tableofcontents

\section{Введение}

Пусть $\Omega$ --- конечное множество мощности $|\Omega|=n$, $G$ --- группа подстановок,
действующая на $\Omega$, т.е.\ подгруппа полной симметрической группы $\Sym(\Omega)$.
Обозначим через $\Set{m}$ множество всех $m$-элементных {\em подмножеств}\ в $\Omega$, а
через $\MSet{m}$ --- множество всех $m$-элементных (или {\em веса $m$}) {\em мультимножеств}\ на $\Omega$,
т.е.\ {\em функций кратности}\ $\varkappa:\,\Omega\to\Zp$ таких, что $\sum_{\omega\in \Omega} \varkappa(\omega)=m$; здесь $\Zp$ --- множество неотрицательных целых чисел. 
Группа  $G$ естественным образом действует на $\Set{m}$ и $\MSet{m}$. Нас интересует оценка числа орбит этих действий при больших $n$ и $m$. Основной тезис работы состоит в том, что во многих естественных ситуациях действие {\em асимптотически свободно}: средний размер стабилизатора множества $x\in \Set{m}$ или мультимножества $x\in \MSet{m}$ близок к единице,
стабилизатор $G_x$ для большинства таких $x$ состоит из тождественного преобразования, а орбиты
$Gx$ {\em регулярны}, т.е.\ имеют максимальный возможный размер, $|G|$.
Соответственно, числа орбит $\Size{m}{G}=|\Set{m}/G|$ и $\MSize{m}{G}=|\MSet{m}/G|$ близки к {\em минимально возможным}, соответствующим свободному действию, а именно, $\Size{m}{G}\sim\eSize{m}/|G|\,$ и $\,\MSize{m}{G}\sim\eMSize{m}/|G|$, где $\eSize{m}=|\Set{m}|$, $\eMSize{m}=|\MSet{m}|$. 

Асимптотическая свобода действия группы перестановок на (мульти)множествах --- явление, давно подмеченное в отдельных содержательных задачах перечислительной комбинаторики. Наиболее ранний известный автору из литературы пример%
\footnote{По свидетельству Ф.~Харари \cite[гл.~9]{Harary-Palmer}, он узнал об этом результате из письма Д.~Пойа в начале 1950-х гг. В статье \cite{FordUhlenbeck1957-4} (p.~167, $2^\circ$) читаем: <<R.J.~Riddell, ``Contributions to the Theory of Condensation'' (dissertation, University of
Michigan, 1951); R.J.~Riddell and G.E.~Uhlenbeck, ``On the Theory of the Virial Development
of the Equation of State of Monoatomic Gases,'' J.\ Chem.\ Phys., 21, 2056, 1953. The contributions
of Polya are given in these references but have not been published by Polya himself.>>}
--- результат об асимптотике 
числа $g_{v,m}$ неориентированных простых графов (с точностью до изоморфизма) c заданными числами вершин $|V|=v$ и ребер $m$.
Здесь $\Omega=V_2$ --- множество пар вершин, $n=|\Omega|={v\choose 2}$. Элементы множества $\Set{m}$ соответствуют помеченным графам с $m$ ребрами; $G=\Sym(V)$ --- группа перенумераций вершин графа; факторизация $\Set{m}$ по действию $G$ соответствует отождествлению графов, отличающихся лишь порядком меток, так что $g_{v,m}=|\Set{m}/G|$. Если отношение  $m/|\Omega|$ --- плотность единиц в матрице инцидентности графа --- не слишком близко к $0$ и $1$, то $g_{v,m}\sim |\Set{m}|/|G|$, т.е.\ действие почти свободно. 
В этом примере и вообще нередко задача перечисления помеченных структур 
несложна, а трудность состоит в количественном учете отождествлений вследствие
<<устранения меток>>. 

Богатый материал для изучения вопроса о существовании и типичности регулярных орбит дает линейная алгебра над конечными полями: действия линейных (аффиных, проективных) групп на конфигурациях
точек, грассманианах и т.д.
В \cite{Hoffman-Welch1968} доказано существование подмножеств линейного пространства $V$ размерности $\geq 5$ над $\FF_2$ с единичным стабилизатором%
\footnote{в оригинале --- {\em totally variant}} 
относительно действия полной линейной группы,
т.е.\ установлен факт существования регулярных орбит $GL(V)$ в $2^V$.  
В работе \cite{Strazdins1997} (p.~166), говоря о классах эквивалентности булевых функций относительно линейных групп, автор замечает (основываясь на фольклоре или интуиции?) что <<асимптотически почти все функции и типы имеют единичный стабилизатор>>%
, не приводя точных формулировок и ссылок.

К линейной алгебре относятся результаты Уальда \cite{Wild00,Wild05} и Хоу \cite{Hou07a,Hou07b,Hou09} об асимптотическом перечислении матроидов и $q$-ичных кодов; показано, что действие группы подстановок элементов фиксированного базиса на грассманиане в $\FF_q^n$ асимптотически ($n\to\infty$) свободно. 
В статье \cite{Lax04} на ту же тему затронута связь с асимптотической теорией характеров ---  привлекательная и не отраженная в настоящей работе точка зрения. 

В ряде работ --- см., например, \cite{Cameron1977,Cameron1999,Siemons1984} --- изучаются другие аспекты действия $G$ на множествах $\Set{m}$, включая и случай бесконечного множества $\Omega$  ($n=\infty$) при конечных числах орбит 
$\Size{m}{G}$. 
В \cite{LW1965} при $n<\infty$ 
доказано, что последовательность $\Size{m}{G}$ унимодальна, а именно: она не убывает при $m\leq n/2$, а при $m>n/2$ имеем $\Size{m}{G}=\Size{n-m}{G}$.

Короткая статья Ж.-М.~Лаборда \cite{Laborde1985} (по заглавию ее нелегко ассоциировать с нашей темой) --- единственная известная автору работа, направленная на выявление общего механизма почти свободного действия. 
Результаты \cite{Laborde1985} изложены в \S~2  п.~4 и в \S~4 п.~4.

Настоящая работа возникла как ответвление в проекте по перечислению тонких градуировок на алгебрах Ли \cite{KPS13}. 
Обсуждение возможности усиления теоремы 6.2 из \cite{KPS13} (публикуется отдельно) привело к задаче о нетривиальной оценке среднего размера стабилизатора.

О структуре работы. 
Работа состоит из <<теоретической части>> (\S\S~2--4) и  приложений (\S~\ref{sec:examples}). 
Итоговые результаты теоретической части --- теоремы~\ref{thm:number_of_orbits} и \ref{thm:number_of_orbits_2tier} об оценках среднего размера стабилизатора (удобным образом нормированного числа орбит) действия
$G$ на $\Set{m}$ и $\MSet{m}$. 
Подсчет числа орбит базируется на формуле Бернсайда,
поэтому наш анализ начинается с оценки числа подмножеств или мультимножеств, инвариантных относительно одной подстановки, в терминах ее цикловой структуры. Этому посвящен \S~2: 
в нем сосредоточена вся <<аналитическая>> часть. Наши неравенства просты, универсальны, но не оптимальны.
В частности, оценка Лаборда числа {\it всех}\ инвариантных подмножеств подстановки
не является <<интегральным>> следствием наших оценок, относящихся к подмножествам {\it фиксированной}\ мощности.
(Но оценки, полученные в \cite{Laborde1985} для последних, слабее наших.)  

В \S~3 результаты, полученные для одной подстановки, распространяются на произвольные совокупности подстановок.
В \S~4 подробно рассмотрен важнейший случай, когда совокупность --- это вся группа $G$. 
Развитие теории в этих параграфах происходит за счет расширения круга задействованных понятий, а доказательства в основном
сводятся к перефразировкам.
Большинство используемых понятий стандартно, но укажем два исключения. (1) В \S~3 определяется цикловой индекс для произвольной совокупности подстановок; в случае групп он отличается от циклового индекса теории перечисления Редфилда-Пойа отсутствием нормализующего множителя. (2) Часто встречающийся объект ---
пара $(\sigma, x)$, где $x$ --- неподвижная точка подстановки $\sigma$, --- кажется, не имеет общепринятого названия. Мы используем термин <<пассивная пара>>. 

Теоремы~\ref{thm:number_of_invariant_sets-1}, \ref{thm:number_of_invariant_sets-K}, \ref{thm:number_of_invariant_sets-mindeg} и \ref{thm:number_of_orbits} составляют логическую
цепочку, в которой происходят последовательные огрубления оценок и, таким образом, потенциально присутствует возможность более полного использования доступной информации о цикловой структуре отдельных подстановок и их совокупностей.  
Вариацию на эту тему представляет ответвление <<Теорема~\ref{thm:number_of_invariant_sets-mindeg} $\Rightarrow$ Следствие~\ref{cor:mindeg-2tier} $\Rightarrow$ Теорема~\ref{thm:number_of_orbits_2tier}>>. 

В дополнение к основным результатам, предложение~\ref{prop:conditions-mn} дает достаточные условия того, что неравенства
теоремы~\ref{thm:number_of_orbits} лучше тривиальных. 
Предложение~\ref{prop:regorbits} --- строгая формулировка импликации <<средний размер стабилизатора близок к единице $\Rightarrow$ почти все орбиты регулярны>>.

Из трех разделов \S~\ref{sec:examples} в двух описаны примеры применения теорем~\ref{thm:number_of_orbits}  и~\ref{thm:number_of_orbits_2tier} к конкретным группам, а в разделе \ref{ssec:gengroups} обсуждаются следствия
теоремы~\ref{thm:number_of_orbits} в свете результатов общей теории групп подстановок, восходящих к работам XIX века Жордана и Бохерта.           
 
Теорема~\ref{thm:GL} в разделе \ref{ssec:affine} о действии аффинной группы $AGL(d,q)$ на точечных конфигурациях в векторном пространстве
$\FF_q^d$ --- иллюстрация эффективности теоремы~\ref{thm:number_of_orbits} в простой ситуации, которая, однако, может представлять достаточно общий интерес. Автору неизвестны публикации, содержащие этот результат. Одно из приложений теоремы~\ref{thm:GL} --- упомянутое выше улучшение оценки из \cite{KPS13}.
Было бы интересно сопоставить теорему~\ref{thm:GL} с техникой и результатами работ \cite{Wild00}--
\cite{Lax04}. 

Цель раздела \ref{ssec:graphs}  --- на примере проиллюстрировать сходство и различие применения теорем~\ref{thm:number_of_orbits} и \ref{thm:number_of_orbits_2tier} (к асимптотическому перечислению, соответственно, однородных гиперграфов и графов).
Теорема~\ref{thm:graphs} значительно перекрывает теорему Обершельпа \cite[\S~9.3]{Harary-Palmer}, но слабее точных результатов, принадлежащих 
А.Д.~Коршунову\cite{Korshunov1970} и Е.М.~Райту \cite{Wright1971,Wright1972}.     

В разделе \ref{ssec:gengroups} мы изучаем, что дают при подстановке в  теорему~\ref{thm:number_of_orbits} известные неравенства, связывающие степень группы подстановок с ее порядком и минимальной степенью.
Неожиданно общий вывод --- теорема~\ref{numorbits_primitive_groups}(б) --- получается для $2$-транзитивных групп, за исключением 
полной симметрической и знакопеременной: такие группы действуют почти свободно при $n=|\Omega|\to\infty$ на $\Set{m}$ и на $\MSet{m}$, во всяком случае, если отношение $m/n$ (сооотв.\ $m/(m+n)$) отделено от $0$ и $1$. 
Известен результат такого же типа \cite{Cameron1986}, относящий к действию на $2^\Omega$ примитивных (на $\Omega$) групп и опирающийся на оценку Лаборда (\ref{Laborde1}).   

\section{Оценки для одной подстановки}
\label{sec:est1}

Оценим число подмножеств и мультимножеств данного размера, инвариантных относительно данной подстановки, в терминах длины циклов этой подстановки. 

\medskip
{\bf 1.} Напомним обозначения: $n=|\Omega|$, $\eSize{m}=|\Set{m}|$ (количество всех $m$-элементных подмножеств множества $\Omega$, $0\leq m\leq n$), и $\eMSize{m}=|\MSet{m}|$ (количество всех $m$-мультимножеств в $\Omega$, $m\geq 0$). Имеем
$$
\eSize{m}=\binom{n}{m},
\qquad
\eMSize{m}=\binom{n+m-1}{m}.
$$ 
Зафиксируем обозначения для отношений $m/n$ и $m/(n+m)$:
\beq{lambda-eta}
 \frac{m}{n}=\lambda,\qquad
 \frac{m}{n+m}=\eta.
\eeq
Введем обозначение для энтропии двузначного распределения с вероятностями 
$p$ и $q=1-p$,
\beq{entropyh}
 h(p)=-p\ln p-q\ln q.
\eeq  
Формула Стирлинга в виде \cite{AbramowitzStegun} (6.1.38)
$$
 n!=\sqrt{2\pi n}\, (n/e)^n e^{\theta_n/(12n)},\quad 0<\theta_n<1,
$$
влечет неравенства, которые понадобятся в доказательстве теоремы \ref{thm:number_of_orbits}:
\begin{equation}
\label{stirling-binomial-1}
 \eSize{m}>\frac{e^{nh(\lambda)}}{\sqrt{2\pi n\lambda(1-\lambda)}}\,e^{-1/6}
\qquad(0<m<n)
\end{equation}
и 
\begin{equation}
\label{stirling-binomial-2}
\eMSize{m}>\frac{e^{(n+m)h(\eta)}}{\sqrt{2\pi (n+m)}}\,
\sqrt{\frac{1-\eta}{\eta}}\,e^{-1/6} \qquad(m>0).
\end{equation}

\medskip
{\bf 2.} Для подстановки $\sigma\in\Sym(\Omega)$ пусть $j_k(\sigma)$ --- число циклов
длины $k$.
Удобный вспомогательный объект --- {\em цикловой моном}  
\beq{zmon}
z_\sigma(t_1,\dots,t_n)=\prod_{k\geq 1} t_k^{j_k(\sigma)}.
\eeq
Здесь $t_1, t_2,\dots, t_n$ --- формальные переменные и 
$\sum_{k=1}^n  k j_k=n$, так что $z_\sigma(t,t^2,\dots,t^n)=t^n$.
Положим
\beq{alpha}
 \alpha(\sigma)=j_1(\sigma)
\eeq
и
\beq{beta}
 \beta(\sigma)=n-\alpha(\sigma)=\sum_{k\geq 2} k j_k(\sigma).
\eeq
Таким образом, $\alpha(\sigma)=|\mathrm{Fix}(\sigma)|$ --- число неподвижных точек подстановки $\sigma$. Число $\beta(\sigma)$ --- размер {\it носителя}\ $\mathrm{supp}(\sigma)=\Omega\setminus\mathrm{Fix}(\sigma)=\{\omega\in\Omega\,|\,\sigma(\omega)\neq\omega\}$ --- будем называть {\em степенью подвижности}%
\footnote{Калька <<степень подстановки>> с английского {\em degree of permutation}\ \cite{Wielandt} двусмысленна; предлагаемый созвучный термин
устраняет двусмысленность, сохраняя
этимологию важного 
понятия {\em минимальной степени}\ совокупности подстановок (\S~\ref{sec:permset}).}
подстановки $\sigma$.

Введем обозначения для числа неподвижных точек действия $\sigma$ на $\Set{m}$ и $\MSet{m}$:  
\beq{alpham}
 \alpha_m(\sigma)=|\Fix{m}(\sigma)|,
 \qquad
 \alpha_{(m)}(\sigma)=|\MFix{m}(\sigma)|,
\eeq
где $\Fix{m}(\sigma)\subset 2^{\Set{m}}$  --- множество $\sigma$-инвариантных $m$-подмножеств в $\Omega$,
а $\MFix{m}(\sigma)\subset 2^{\MSet{m}}$ ---  множество $\sigma$-инвариантных $m$-мультимножеств. 

Производящие функции для величин (\ref{alpham}) суть
\beq{Asigma}
 A_\sigma(t)=\sum_{m=0}^n \alpha_m(\sigma) t^m
\eeq
и
$$
 A_{(\sigma)}(t)=\sum_{m=0}^\infty \alpha_{(m)}(\sigma) t^m.
$$
Хорошо известно и элементарно выводится тождество
$$
 A_\sigma(t)=\prod_{k\geq 1} (1+t^k)^{j_k(\sigma)}= z_\sigma(1+t,1+t^2,1+t^3,\dots).
$$
В самом деле, $\Omega$ есть объединение $\sigma$-орбит (циклов подстановки $\sigma$). 
Всякое $\sigma$-инвариантное подмножество $x\subset\Omega$ однозначно соответствует подмножеству
$\bar x$ пространства орбит $\Omega/\sigma$. 
При раскрытии скобок в произведении $P=\prod_{\xi\in \Omega/\sigma} (1+\tau_\xi)$, где $\tau_\xi$ --- формальные переменные,  подмножеству $\bar x$ соответствует единственный моном $p_x=\prod_{\xi\in \bar x} \tau_\xi$. Подставляя $\tau_\xi=t^{|\xi|}$, получим $P=\prod_{k\geq 1} (1+t^k)^{j_k(\sigma)}$,
а моном $p_x$ обращается в $t^m$ с показателем $m=\sum_{\xi\in\bar x}|\xi|=|x|$, равным размеру подмножества $x$.   

Аналогично имеет место тождество 
$$
 A_{(\sigma)}(t)=\prod_{k\geq 1} (1-t^k)^{-j_k(\sigma)}= z_\sigma\left(\frac{1}{1-t},\frac{1}{1-t^2},\dots\right).
$$
В частности, когда $\sigma=e$ (тождественная подстановка), имеем
$$
 A_e(t)=(1+t)^n,\qquad
 A_{(e)}(t)=(1-t)^{-n}.
$$

\medskip
{\bf 3.} Отправной точкой для последующего анализа являются оценки числа $\sigma$-инвариантных $m$-(мульти)множеств. Используемыe в теореме~\ref{thm:number_of_invariant_sets-1} обозначения введены формулами
(\ref{lambda-eta})--(\ref{entropyh}) и  
(\ref{beta})--(\ref{alpham}). 
Положим также
\beq{xi12}
\begin{array}{l}\displaystyle
\xi_1=\xi_1(\lambda)=\left(1-2\lambda(1-\lambda)\right)^{1/2}=\frac{\sqrt{m^2+(n-m)^2}}{n}, 
\\[2ex]\displaystyle
\xi_2=\xi_2(\eta)=\left(\frac{1-\eta}{1+\eta}\right)^{1/2}=\left(\frac{n}{n+2m}\right)^{1/2}.
\end{array}
\eeq

\begin{theorem}
\label{thm:number_of_invariant_sets-1}
Для подстановки $\sigma$ справедлива верхняя оценка числа $\sigma$-инвариантных $m$-подмножеств 
	\beq{ub_alpham}
 \alpha_m(\sigma) 
\,\leq\,
 \xi_1^{\beta(\sigma)}
\,e^{n h(\lambda)}
	\eeq
и верхняя оценка числа $\sigma$-инвариантных $m$-мультимножеств
\beq{ub_malpham}
 \alpha_{(m)}(\sigma) 
\,\leq\,
\xi_2^{\beta(\sigma)}
\,e^{(n+m) h(\eta)}.
\eeq
\end{theorem}

\proof
(а) При любом $\tau>0$ функция $ k\;\mapsto\;(1+\tau^k)^{1/k}$ убывает.
Следовательно, справедливо неравенство
\beq{power_inequality-1}
 1+\tau^k\leq (1+\tau^2)^{k/2},\quad k\geq 2.
\eeq
Применяя это неравенство с $\tau=t^3,t^4,\dots$, находим
$$
 A_\sigma(t)\leq 
(1+t)^{\alpha(\sigma)} (1+t^2)^{\beta(\sigma)/2}=A_e(t)\left(\frac{1+t^2}{(1+t)^2}\right)^{\beta(\sigma)/2}.
$$
Поскольку коэффициенты производящей функции $A_\sigma(t)$ неотрицательны, имеем:
для любого $t>0$
\beq{alpha_m_estimate}
\alpha_m(\sigma) \leq t^{-m} A_\sigma(t) \leq t^{-m}A_e(t)\,\left(\frac{1+t^2}{(1+t)^2}\right)^{\beta(\sigma)/2}.
\eeq
Минимум функции $t^{-m} A_e(t)=t^{-m} (1+t)^n$ достигается в точке
$$
 t_*=\frac{m}{n-m}=\frac{\lambda}{1-\lambda}. 
$$   
Вычисляя, находим
$$
\frac{1}{1+t_*}=1-\lambda,\qquad
\frac{t_*}{1+t_*}=\lambda
$$
и 
$$
 t_*^{-m} A_e(t_*)=\left[\lambda^{\lambda}(1-\lambda)^{1-\lambda}\right]^{-n} = e^{n h(\lambda)}.
$$ 
Имеем также
$$
 \frac{1+t_*^2}{(1+t_*)^2}=1-\frac{2t_*}{(1+t_*)^2}
=1-2\lambda(1-\lambda)=\xi_1^2. 
$$
Используя значение $t=t_*$ в неравенстве (\ref{alpha_m_estimate}), придем к оценке
(\ref{ub_alpham}).

\medskip
(б) Аналогично, для вывода оценки (\ref{ub_malpham}) используем неравенство
$$
(1-\tau^k)^{-1}\leq (1-\tau^2)^{-k/2},\quad k\geq 2,
$$
справедливое при любом $\tau\in [0,1)$.
Подставляя $\tau=t^3,t^4,\dots$, находим
$$
 A_{(\sigma)}(t)\leq 
(1-t)^{\alpha(\sigma)} (1-t^2)^{-\beta(\sigma)/2}=A_{(e)}(t)\left(\frac{1-t}{1+t}\right)^{\beta(\sigma)/2}.
$$
Как и в (\ref{alpha_m_estimate}), имеем: при любом $t\in (0,1)$ 
\beq{malpha_m_estimate}
\alpha_{(m)}(\sigma) \leq t^{-m} A_{(\sigma)}(t) \leq t^{-m}A_{(e)}(t)\,\left(\frac{1-t}{1+t}\right)^{\beta(\sigma)/2}.
\eeq
В правой части в качестве $t$ используем значение 
$$
 t_{(*)}=\frac{m}{n+m}=\eta,
$$   
доставляющее минимум функции $t^{-m} A_{(e)}(t)=t^{-m} (1-t)^{-n}$.
Значение минимума есть $t_{(*)}^{-m} A_{(e)}(t_{(*)})=e^{(m+n)h(\eta)}$,
и теперь из (\ref{malpha_m_estimate}) следует (\ref{ub_malpham}).
\eop

\bigskip 
{\bf 4.} 
Для полного числа $\sigma$-инвариантных подмножеств имеет место оценка Ж.-М.~Лаборда
\cite[теорема~1]{Laborde1985} 
\beq{Laborde1}
\sum_{m=1}^n \alpha_m(\sigma)\leq 2^{n-\beta(\sigma)/2}.
\eeq
Она не является следствием  
теоремы \ref{thm:number_of_invariant_sets-1},
но 
доказывается проще: используя производящую функцию
(\ref{Asigma}), имеем 
$$
 \sum_{m=1}^n \alpha_m(\sigma)=A_\sigma(1)=z_\sigma(2,2,\dots,2)=2^n\cdot 2^{-\sum_{k=2}^n (k-1)j_k},
$$
откуда ввиду (\ref{beta}) получаем (\ref{Laborde1}).

\section{Оценки для совокупности подстановок}
\label{sec:permset}

Назовем пару $(\sigma, x)$, где $x$ --- мультимножество (в частности, подмножество) на $\Omega$,  
{\em пассивной парой},
если $\sigma(x)=x$. 
Наша цель в этом разделе --- получить общие оценки для числа пассивных пар $(\sigma, x)$ с определенными ограничениями на $\sigma$ и $x$.

Пусть $K$ --- подмножество (не обязательно подгруппа) группы $\Sym(\Omega)$.
Вся информация о цикловой структуре подстановок совокупности $K$ содержится в {\em цикловом индексе} --- сумме цикловых мономов вида (\ref{zmon}) 
$$
 z^K(t_1,\dots, t_n)=\sum_{\sigma\in K} z_\sigma(t_1,\dots,t_n).
$$
Нам понадобится производящая функция более частного вида, содержащая лишь информацию о распределении числа неподвижных точек или, эквивалентно, о распределении размеров носителей. 
Пусть $f^K_{k}$ --- количество таких $\sigma\in K$, для которых $\beta(\sigma)=k$.
Положим
$$
 F^K(t)=\sum_{k=0}^n f^K_k t^k.
$$
Легко видеть, что
$$
 F^K(t)=\sum_{\sigma\in K} t^{\beta(\sigma)} = z^K(1,t^2,t^3,\dots).
$$

Пусть $u^K_m$ (соотв.,\ $u^K_{(m)}$) --- число таких пассивных пар $(\sigma, x)$, для которых $\sigma\in K$ и $x\in\Set{m}$ (соотв.,\ $\sigma\in K$, $x\in\MSet{m}$).
Очевидны равенства
$$
 u^K_m=\sum_{\sigma\in K} \alpha_m(\sigma),
\qquad
 u^K_{(m)}=\sum_{\sigma\in K} \alpha_{(m)}(\sigma).
$$

%
Сумма по $\sigma\in K$ оценок из теоремы~\ref{thm:number_of_invariant_sets-1} дает следующий результат.

\begin{theorem}    
\label{thm:number_of_invariant_sets-K}
 Имеют место неравенства
\begin{equation}
\label{num_inv_pairs-1}
u^K_m\leq e^{n h(\lambda)}\, F^K(\xi_1(\lambda))
\end{equation}
и
\begin{equation}
\label{num_inv_pairs-2}
u^K_{(m)}\leq e^{(n+m) h(\eta)}\,  F^K(\xi_2(\eta)).
\end{equation}
\end{theorem}

\smallskip\noindent{\bf Определение. }
Назовем {\em минимальной степенью подмножества}\ $K\subset\Sym(\Omega)$ минимальный размер носителя подстановки
$\sigma\in K$:
$$
 \mu^K=\min_{\sigma\in K} \beta(\sigma).
$$

\smallskip\noindent{\bf Замечание. }
Если $e\in K$, наше определение дает $\mu^K=0$. В частности, так обстоит дело, когда $K$ --- группа. 
В этом случае для разрешения терминологической коллизии со стандартным определением (см., напр., \cite{Wielandt}) приходится разделять понятия {минимальной степени подмножества}\ и {минимальной степени группы}.  
{\em Минимальная степень группы}\ $G$ в смысле стандартного определения равна минимальной степени подмножества $\Gminuse=G\setminus\{e\}$ в нашем смысле.
В явном виде, 
\begin{equation}
\label{mindeg_group}
\mu^{\Gminuse}=\min_{g\in G, g\neq e}\big|\{\omega\in\Omega\,|\,g(\omega)\neq\omega\}\big|.
\end{equation}

\begin{theorem}
\label{thm:number_of_invariant_sets-mindeg}
Числа пассивных пар оцениваются сверху через минимальную степень $\mu=\mu^K$ подмножества $K$
следующим образом:
\begin{equation}
\label{num_inv_pairs-mindeg1}
u^K_m\leq |K|\,e^{n h(\lambda)}\, \xi_1^{\mu}
\end{equation}
и
\begin{equation}
\label{num_inv_pairs-mindeg2}
u^K_{(m)}\leq |K|\,e^{(n+m) h(\eta)}\, \xi_2^{\mu}.
\end{equation}
\end{theorem}

\proof
Используем теорему~\ref{thm:number_of_invariant_sets-K} и оценку $F^K(\xi)\leq |K|\, \xi^{\mu}$. 
\eop

\medskip
Величины $u^K_m$ и $u^K_{(m)}$, рассматриваемые как функции множества $K$, аддитивны. Это очевидное замечание может быть полезно для изоляции влияния <<плохого>>
множества $B$ малого размера на оценку в целом:

\begin{Corollary}
\label{cor:mindeg-2tier}
В условиях теоремы~{\rm\ref{thm:number_of_invariant_sets-mindeg}} пусть
$B\subsetneq K$ и $\mu^*:=\mu^{K\setminus B}>\mu$. Тогда оценки {\rm (\ref{num_inv_pairs-mindeg1})} 
и {\rm (\ref{num_inv_pairs-mindeg2})} можно улучшить:
\begin{equation}
\label{num_inv_pairs-twotier1}
u^K_m\leq 
e^{n h(\lambda)}\,\left( (|K|-|B|)\,\xi_1^{\mu^*} 
+|B|\,\xi_1^{\mu}\right)
\end{equation}
и
\begin{equation}
\label{num_inv_pairs-twotier2}
u^K_{(m)}\leq 
e^{(n+m) h(\eta)}\,\left( (|K|-|B|)\,\xi_2^{\mu^*} 
+|B|\,\xi_2^{\mu}\right).
\end{equation}
\end{Corollary}

\section{Оценки среднего размера стабилизатора}
\label{sec:avstab}

Здесь рассмотрим оценки, полученные в \S~\ref{sec:permset}, с новой точки зрения в случае, когда подмножество $K\subset\Sym(\Omega)$ является группой, $K=G$.

\medskip
{\bf 1.} Пусть $S$ --- конечное множество, на котором действует группа $G$,
и $S/G$ --- множество орбит. Для числа орбит $|S/G|$ имеем цепочку равенств
$$
 |S/G|=\dst
\sum_{o\in S/G} 1=\frac{1}{|G|}\sum_{o\in S/G}
\sum_{s\in o} |\Stab(s)|
\;\stackrel{*}{=}\;\frac{1}{|G|}\dst\sum_{\sigma\in G} |\Fix{S}(\sigma)|,
$$
где $\Fix{S}(\sigma)=\{s\in S\,|\,\sigma(s)=s\}$. 
Шаг $(*)$ --- это подсчет числа пассивных пар
$(\sigma,s)\in G\times S$ двумя способами, а результат --- 
известная формула Бернсайда (или Коши-Фробениуса-Бернсайда).

Вводя {\em средний размер стабилизатора}\ элемента $s\in S$
$$
 \avStab{G}{S}=\frac{1}{|S|}\sum_{s\in S} |\Stab(s)|,
$$
можем написать
$$
 |S/G|=\frac{|S|}{|G|}\avStab{G}{S}.
$$ 
Очевидным образом $1\leq\avStab{G}{S}\leq |G|$; соответственно, 
\beq{triv_num_orbits}
\frac{|S|}{|G|}\leq |S/G|\leq |S|.
\eeq
Верхняя граница отвечает случаю {\it тривиального}\ действия, когда все орбиты одноточечны. Нижняя 
граница отвечает случаю {\it свободного}\ действия, когда все подстановки $\sigma\in G$ действуют без неподвижных точек. 

В случае $|S|\leq|G|$ тавтологичное неравенство $|S/G|\geq 1$ (равенство, если действие транзитивно) влечет $\avStab{G}{S}\geq |G|/|S|$.   

\medskip
{\bf 2.}
Пусть теперь множество $S$ есть $\Set{m}$ или $\MSet{m}$.
Объектом нашего внимания будут улучшения верхней оценки в (\ref{triv_num_orbits}), в частности --- условия, при которых 
средний размер стабилизатора близок к единице.
Считая группу $G$ и множество $\Omega$ фиксированными, будем писать $\avStabSet{m}$
вместо $\avStab{G}{\Set{m}}$ и $\avStabMSet{m}$
вместо $\avStab{G}{\MSet{m}}$. 

\begin{theorem}    
\label{thm:number_of_orbits}
Пусть $\mu=\mu^{\Gminuse}$ --- минимальная степень группы подстановок $G$, действующей на множестве $\Omega$ мощности $n$. 
Для естественного действия группы $G$ на подмножествах мощности $m\neq 0,n$ в $\Omega$ имеет место оценка
среднего размера стабилизатора
\begin{equation}
\label{general_stab_size-sets}
\frac{\avStab{}{m}-1}{|G|}
<
3\sqrt{\frac{m(n-m)}{n}}\,\cdot\xi_1^{\mu}. 
\end{equation}
Аналогично, для естественного действия группы $G$ на мультимножествах веса $m>0$
имеет место оценка
\begin{equation}
\label{general_stab_size-msets}
\frac{\avStab{}{(m)}-1}{|G|}
<
3\sqrt{\frac{m(n+m)}{n}}\,\cdot\xi_2^{\mu}.
\end{equation}
В этих оценках $\xi_1$ и
$\xi_2$ --- числа между $0$ и $1$, определенные в {\rm (\ref{xi12})}.
\end{theorem}

\proof
В обозначениях \S\S~\ref{sec:est1}--\ref{sec:permset}  ($\eSize{m}=|\Set{m}|$, $u^G_m$ --- число пассивных пар в
$G\times \Set{m}$, $\Gminuse=G\setminus\{e\}$ и т.д.)
средний размер стабилизатора и его отклонение от единицы в двух интересующих нас случаях суть
$$
\ba{ll}
\dst
\avStab{}{m}=\frac{u^G_m}{\eSize{m}}, 
\quad &\dst
\avStab{}{m}-1=\frac{u^{\Gminuse}_m}{\eSize{m}}, 
\\[2ex]
\dst
\avStab{}{(m)}=\frac{u^G_{(m)}}{\eMSize{m}}, 
\quad &\dst
\avStab{}{(m)}-1=\frac{u^{\Gminuse}_m}{\eMSize{m}}.
\ea
$$
Остается применить теорему~\ref{thm:number_of_invariant_sets-mindeg} с учетом неравенств
(\ref{stirling-binomial-1}), (\ref{stirling-binomial-2}), тождеств
$n\lambda(1-\lambda)=m(n-m)/n$, $\eta/(1-\eta)=m/n$ и числового неравенства
$e^{1/6}\sqrt{2\pi}< 3$.
\eop

\bigskip
{\bf 3.}
Насколько информативна теорема~\ref{thm:number_of_orbits}?
В силу известных результатов теории групп подстановок --- см.\ неравенства  (\ref{Babai_ineq_mindeg}) и (\ref{Bochert_ineq_mindeg}) в \S~\ref{ssec:gengroups} ---  для транзитивных групп, кроме $Sym(\Omega)$ и $Alt(\Omega)$ при всех $n\geq 10$ верно соотношение $\mu>\frac{1}{2}\sqrt{n}$. 
С учетом этого приведем достаточные асимптотические условия того, что правые части неравенств (\ref{general_stab_size-sets}) и (\ref{general_stab_size-msets}) 
произвольно малы. 
Ниже {\it почти все}\ означает <<все, за исключением, возможно, конечного числа>>.

\begin{proposition}
\label{prop:conditions-mn} 
Для заданного $c>0$ пусть 
$$
\ba{l}
 R_c=\{(m,n)\in \ZZ_+^2\,|\, m>c\sqrt{n}\ln n\},
\\[1ex]
\hat R_c=\{(m,n)\in \ZZ_+^2\,|\, \min(m,n-m)>c\sqrt{n}\ln n\}.
\ea
$$
Если $c>1/2$, то для любого $\eps>0$ 

\smallskip\noindent
{\rm(а)} 
при почти всех $(m,n)\in\hat R_c$ выполнено неравенство
\beq{suffcond-applicability-sets}
\sqrt{\frac{m(n-m)}{n}}\,\cdot\left(\xi_1\left({m}/{n}\right)\right)^{\sqrt{n}/2}<\eps;
\eeq
%
{\rm(б)} 
при почти всех $(m,n)\in R_c$ таких, что $n>16$, выполнено неравенство
\beq{suffcond-applicability-msets}
\sqrt{\frac{m(n+m)}{n}}\,\cdot\left(\xi_2\left(\frac{m}{m+n}\right)\right)^{\sqrt{n}/2}<\eps.
\eeq
\end{proposition}

\proof
(а) Пусть для определенности $m\leq n/2$. 
Ввиду неравенства 
$\ln\xi_1(\lambda)<-\lambda(1-\lambda)$
логарифм левой части неравенства (\ref{suffcond-applicability-sets}) оценивается сверху величиной
$$
 A= \frac{1}{2}\left(\ln m-\sqrt{n}\cdot \lambda(1-\lambda)\right)= \frac{\ln m}{2}\left(1-\frac{m/\ln m}{\sqrt{n}}(1-\lambda)\right).
$$ 
Пусть $R_c'$ --- часть множества $R_c$, выделенная условием $m/n<d/2$, где $d=1-(2c)^{-1}>0$, 
а $R_c''$ --- часть множества $R_c$, выделенная условием $d/2\leq m/n\leq 1/2$.
Ясно, что для почти всех $(m,n)\in R_c''$ имеем $A<\ln \eps$.   
Если же 
$(m,n)\in R_c'$, то 
$$
\frac{1-\lambda}{\sqrt{n}}\,\cdot\,\frac{m}{\ln m}\;>\;
\frac{1-d/2}{\sqrt{n}}\,\cdot\,\frac{c\sqrt{n}\ln n }{\frac{1}{2}\ln n+\ln (c\ln n)}.
$$
Предел правой части при $n\to\infty $ равен $(1-d/2)\cdot 2c=c+1/2>1$, поэтому $A<\ln\eps$ почти всюду в $R_c'$.

\medskip
%
(б) Так же, как в п.~(а), доказывается, что неравенство
(\ref{suffcond-applicability-msets})  справедливо при почти всех $(m,n)\in R_c$, для которых $m/n\leq 1/2$.
Новым по сравнению с (а) является случай подмножества $R_c'''\subset R_c$, выделенного условиями $2m>n>16$. 
Поскольку $(m(m+n)/n)^{1/2}\leq (1-\eta)^{-1} n^{1/2}$, левая  часть неравенства (\ref{suffcond-applicability-msets}) оценивается сверху величиной
$$
B=\sqrt{n}\left(\frac{1-\eta}{1+\eta}\right)^{-1+\sqrt{n}/4}.
$$
Условие $2m>n$ влечет $\eta>1/3$, поэтому неравенство $B<\eps$ в $R_c'''$ может нарушаться лишь для конечного числа значений $n$. 
Если $n>16$ --- такое значение, то $B\to 0$ при $\eta\to 1$, значит, неравенство $B\geq \eps$ возможно лишь для конечного числа значений $m$.
Доказательство  завершено.
\eop

\medskip
{\bf 4.}
Оценку (\ref{general_stab_size-sets}) можно сравнить
с теоремой 2 в статье 
\cite{Laborde1985}, доказанной при условии  
$\mu\ge \mu_{\min}(n)\sim n/3$.
В наших обозначениях этот результат выглядит так: 
$$
 \frac{\avStab{}{m}-1}{|G|}<\frac{c_m}{\eSize{m}},
$$
где  $\sum_{m=0}^{n}c_m x^m=(1+x)^p(1+x^2)^{(n-p)/2}$,
$p=\lceil \frac{4}{3}n-{\mu}\rceil$. Показано, что $c_m\leq 2^{p/2-1} \sqrt{\eSize{m}}$. 
Нетрудно видеть, что неравенство  (\ref{general_stab_size-sets}) сильнее.

\medskip
{\bf 5.}
Отметим, что присутствие множителей с квадратным корнем в правых частях оценок 
теоремы~\ref{thm:number_of_orbits}\ --- результат оценки коэффициентов производящих функций $A_{\sigma}(t)$ и $A_{(\sigma)}(t)$ в теореме~\ref{thm:number_of_invariant_sets-1} посредством неравенства Коши. 
Более аккуратные оценки с помощью метода перевала позволили бы (во всяком случае, в конкретных случаях) избавиться от этих множителей.
Более принципиален вопрос о сравнительной величине $\mu$ и $\log |G|$. Грубо говоря, теорема~\ref{thm:number_of_orbits} показывает, что действие $G$ на множествах и мультимножествах подходящего размера почти свободно при условии $\mu\gg\log |G|$. Это условие --- само по себе не экзотическое (оно выполнено для 
2-транзитивных групп подстановок большого размера, см.~(\ref{Bochert_ineq_mindeg}) и (\ref{Pyber_ineq_order})) ---
можно ослаблять за счет использования более детальной информации о циклах подстановок группы $G$. 

Улучшение возможно, например, в ситуации, когда лишь малое число элементов группы (скажем, элементы <<плохого>> подмножества $B$) имеют близкую к $\mu$ или, во всяком случае, недостаточно хорошо оцениваемую снизу степень подвижности, в то время как степень подвижности большинства элементов намного превышает $\mu$, 
то есть $\mu^{G\setminus B}\gg\mu$. Следующая теорема относится к теореме~\ref{thm:number_of_orbits}\ так же, как к теореме~\ref{thm:number_of_invariant_sets-mindeg}
относится следствие~\ref{cor:mindeg-2tier}. 
  
\begin{theorem}    
\label{thm:number_of_orbits_2tier}
Пусть, в дополнение к условиям теоремы~{\rm \ref{thm:number_of_orbits}},
для подмножества $B$ группы $G$ известно, что $\mu^*:=\mu^{G\setminus B}>\mu$.
Тогда 
\begin{equation}
\label{general_stab_size-sets2}
\avStab{}{m}-1
<
3\sqrt{\frac{m(n-m)}{n}}\,
\left(|G|\,\xi_1^{\mu^*}\, 
+\,|B|\,\xi_1^{\mu}\right)
\end{equation}
 и
\begin{equation}
\label{general_stab_size-msets2}
\avStab{}{(m)}-1
<
3\sqrt{\frac{m(n+m)}{n}}
\left(|G|\,\xi_2^{\mu^*}\,+\,|B|\,\xi_2^{\mu}\right).
\end{equation}
Более общим образом, если имеем цепочку вложенных множеств $B_1\subset B_2\subset\dots\subset B_p=G$
и $\mu_i=\mu^{G\setminus B_i}$ $(i<p)$, то
\begin{equation}
\label{general_stab_size-sets3}
\avStab{}{m}-1
<
3\sqrt{\frac{m(n-m)}{n}}\,
\left(\sum_{i=1}^{p-1} (|B_{i+1}|-|B_{i}|)\xi_1^{\mu_i}\,+|B_1|\,\xi_1^\mu\right),
\end{equation}
и аналогичная оценка имеет место для $\avStab{}{(m)}-1$.
\eop
\end{theorem}

{\it Расстояние Хемминга}\ на группе подстановок $G$ с  пространством действия $\Omega$ определяется формулой
$$
 d(\sigma_1,\sigma_2)=\big|\{\omega\in \Omega\,|\,\sigma_1(\omega)\neq \sigma_2(\omega)\}\big|.
$$ 
В метрике Хемминга шар радиуса $r$ с центром в единице группы --- это множество $B(r)$ элементов, степень подвижности которых не превышает $r$.
Сфера радиуса $r$ с центром в единице, $\partial B(r)=B(r)\setminus B(r-1)$, --- это множество элементов, степень подвижности которых в точности равна $r$. 
Используя эти обозначения, получаем элегантное

\begin{Corollary}
Имеют место неравенства
$$
\avStab{}{m}-1
<
3\sqrt{\frac{m(n-m)}{n}}\,
\sum_{r=1}^{n} |\partial B(r)|\xi_1^{r}
$$
и
$$
\avStab{}{(m)}-1
<
3\sqrt{\frac{m(n+m)}{n}}\,
\sum_{r=1}^{n} |\partial B(r)|\xi_2^{r}.
$$
В действительности суммирование начинается с $r=\mu$ (минимальная степень группы $G$), поскольку $\partial B(r)=\emptyset$ при $1\leq r<\mu$. 
\end{Corollary}

\medskip
Успех применения теоремы~{\rm \ref{thm:number_of_orbits_2tier}} в виде (\ref{general_stab_size-sets2}), (\ref{general_stab_size-msets2}) для доказательства асимптотической свободы действия $G$ на (мульти)множествах  зависит от возможности выделить <<плохое, но малое>> подмножество $B\subset G$ такое, что 
$$
\mu^*\gg \log|G|, \qquad \mu\gg \log|B|.
$$ 
Общий вариант теоремы предоставляет б\'ольшую гибкость. Пример его использования --- доказательство теоремы~\ref{thm:graphs} ниже.  

\medskip
{\bf 6.}
Наконец, 
приведем не зависящие от изложенной теории, но
полезные для ее приложений нижние оценки числа {регулярных}\ орбит  и принадлежащих им точек.
Их удобно представить в вероятностных терминах.
(Часть (а) --- это лемма 1 в \cite{Cameron1986}; наше $\delta$ там обозначено $\eps_1(G)$.)   

\begin{proposition}
\label{prop:regorbits} 
Пусть группа $G$ действует на множестве $S$, при этом $\delta$ --- отклонение среднего размера стабилизатора от $1$, 
$$
\delta=\avStab{G}{S}-1.
$$
{\rm(а)} Пусть на пространстве орбит $S/G$ задано равномерное вероятностное распределение (т.е.\ все
орбиты равновероятны). Тогда вероятность того, что случайно выбранная орбита регулярна, не меньше, чем $(1-\delta)/(1+\delta)$.
\\[1ex]
{\rm(б)} Пусть на множестве $S$ задано равномерное вероятностное распределение (т.е.\ все
точки равновероятны). Тогда вероятность того, что случайно выбранная точка $s$ принадлежит 
регулярной орбите (эквивалентно --- что $|\Stab(s)|=1$), не меньше, чем $1-\delta$.
\end{proposition}

\proof
(а) Пусть $R$ --- множество регулярных орбит. Тогда 
$$
o\in R \;\Rightarrow \; |o|=|G|,
\qquad
o\notin R \;\Rightarrow \; |o|\leq|G|/2.
$$
Следовательно,
\beq{regorb}
 |S|=\sum_{o\in S/G} |o| =\sum_{o\in R}+\sum_{o\notin R}
 \leq |R|\,|G|+(|S/G|-|R|)\frac{|G|}{2}.
\eeq
Деля обе части на $|G|$, заменяя слева $|S|/|G|=|S/G|(1+\delta)^{-1}$, получим
$$
 |S/G|\left(\frac{1}{1+\delta}-\frac{1}{2}\right)\leq \frac{|R|}{2}.
$$ 
Значит, доля регулярных орбит $|R|/|S/G|\geq (1-\delta)/(1+\delta)$.

\medskip
(б) Делая в правой части (\ref{regorb}) замену $|S/G|=(1+\delta)|S|/|G|$,
находим
$$
 |R|\,|G|\geq |S|(1-\delta).
$$
Левая часть ---  
	число точек, принадлежащих регулярным орбитам.
\eop

\section{Приложения} 
\label{sec:examples}


\subsection{Число классов аффинной эквивалентности наборов векторов}
\label{ssec:affine}
Возьмем в качестве $\Omega$ линейное пространство $V$ размерности $d$ над конечным полем
$\FF_q$. Перечисление классов аффинной эквивалентности $m$-точечных непомеченных конфигураций в пространстве $V$ равносильно перечислению орбит действия аффинной группы $G=AGL(d,q)$ на множестве $\Set{m}$
(если все точки конфигурации различны), либо на $\MSet{m}$ (если допускаются кратные точки). Следуя обозначениям в начале статьи, назовем эти числа орбит соответственно $\Size{m}{AGL(d,q)}$ и $\MSize{m}{AGL(d,q)}$. 

\smallskip
Пусть $q$ фиксировано, а $d,m\to\infty$. 
Отметим, что
$$
\begin{array}{c}
 n=|\Omega|=q^d,
 \\[1ex]\displaystyle
|G|=|V|\cdot|GL(d,q)|=q^d\prod_{j=0}^{d-1}(q^d-q^j) < q^{d^2+d}. 
\end{array} 
$$

\bigskip
\begin{theorem}
\label{thm:GL}
{\rm (а)} Отклонение от $1$ среднего значения стабилизатора действия группы $AGL(d,q)$ на множестве $m$-точечных непомеченных конфигураций 
без кратных точек $(0<m<n)$ в аффинном пространстве $\FF_q^d$ (мощности $n=q^d$) удовлетворяет оценке
\begin{equation}
\label{linear_stab_size-sets}
\avStab{}{m}-1<
	3\sqrt{\frac{m(n-m)}{n}}\,\cdot\exp \left\{n(1-q^{-1})\ln\xi_1(\lambda)+d^2\ln q\right\}.
\end{equation}
Для действия $AGL(d,q)$ на множестве $m$-точечных непомеченных конфигураций с допустимыми кратными точками $(m>0)$
справедлива оценка
\begin{equation}
\label{linear_stab_size-msets}
\avStab{}{(m)}-1<
	3\sqrt{\frac{m(n+m)}{n}}\,\cdot\exp \left\{n(1-q^{-1})\ln\xi_2(\eta)+d^2\ln q\right\}.
\end{equation}
В этих неравенствах $\lambda=m/n$, $\eta=m/(m+n)$, $\xi_1(\lambda)$ и $\xi_2(\eta)$ --- числа между 0 и 1, определенные формулой {\rm(\ref{xi12})}.

{\rm (б)} Положим $\bar{m}=\min (m, n-m)$.
Если $q$ фиксировано, а $\bar{m}$ и $d$  стремятся к бесконечности и при этом 
$\;\underline\lim (\bar{m}/d^2)> (1-q^{-1})^{-1} \ln q$, то
$$
 \Size{m}{AGL(d,q)}= \frac{\displaystyle{n\choose m}}{|AGL(d,q)|}\,(1+o(1)).
$$
Аналогично, если $q$ фиксировано, а ${m}$ и $d$  стремятся к бесконечности и при этом 
$\;\underline\lim ({m}/d^2)> (1-q^{-1})^{-1} \ln q$, то
$$
 \MSize{m}{AGL(d,q)}= \frac{\displaystyle{n+m-1\choose m}}{|AGL(d,q)|}\,(1+o(1)).
$$
\end{theorem}

\proof
(а) Минимальная степень группы $G$ есть $\mu=q^{d}-q^{d-1}=n(1-q^{-1})$. (Ср.\ \cite{DixonMortimer1996}, Exercise 3.3.3).
Действительно, неподвижные точки аффинного отображения $v\mapsto Av+t$ образуют аффинное подпространство, параллельное линейному подпространству $\mathrm{Ker}(A-I)$; его мощность $\leq q^{d-1}$ при $A\neq I$. Если же $A=I$, но $t\neq 0$, то неподвижных точек нет.

Применение теоремы~\ref{thm:number_of_orbits} ведет к неравенствам (\ref{linear_stab_size-sets}), (\ref{linear_stab_size-msets}).

\medskip
(б) Подлежащая доказательству асимптотика величины $\Size{m}{AGL(d,q)}$ эквивалентна утверждению $\avStab{}{m}-1\to 0$.
Обозначим через $A$ логарифм правой части  (\ref{linear_stab_size-sets}). 
В предположении $\lambda=m/q^d\to 0$ имеем асимптотику
$$
A\sim -m(1-q^{-1})(1+o(1)) + d^2\ln q \to -\infty,
$$ 
согласно условию на $m/d^2$.
Чтобы принять во внимание возможность $\overline{\lim}\,\lambda>0$, достаточно повторить рассуждение, использованное   
в доказательстве части (а) предложения~\ref{prop:conditions-mn}. 

Аналогично доказывается утверждение об асимптотике  $\MSize{m}{AGL(d,q)}$. 
Нижняя грань отношения $m/d^2$ играет ключевую роль при рассмотрении 
случая $\eta=m/(m+q^d)\to 0$; при этом асимптотика логарифма правой части  (\ref{linear_stab_size-msets})  
такая же, как у величины $A$ выше. В остальном копируем рассуждение
из доказательства части (б) предложения~\ref{prop:conditions-mn}.
\eop

\medskip\noindent
{\bf Замечание. } (а) При переходе от группы $G$ к ее подгруппе, действующей на том же множестве, размеры стабилизаторов не увеличиваются. Следовательно, {\em оценки {\rm(\ref{linear_stab_size-sets})}
и {\rm(\ref{linear_stab_size-msets})} справедливы для любой подгруппы полной аффинной группы, в частности, для любой линейной группы, действующей на пространстве $V=\FF_q^d$. }
Назовем теоремой~$\ref{thm:GL}'$ копию теоремы~\ref{thm:GL} с группой $G=GL(d,q)$ вместо $G=AGL(d,q)$.

\smallskip
(б) Напротив, при ограничении действия группы $G$ на инвариантное подмножество $\Omega'\subset\Omega$ нет оснований для того, чтобы имеющиеся для $\Set{m}$ и $\MSet{m}$ оценки среднего размера стабилизаторов перенести на $\Omega'_{m}$ и $\Omega'_{(m)}$. 

\smallskip
Например, не проходит следующая естественная попытка вывести теорему~\ref{thm:GL} из
теоремы~$\ref{thm:GL}'$. Выбрав не содержащую $0$ гиперплоскость $V$ в линейном пространстве $\Omega=\FF_q^{d+1}$, 
реализуем полную аффинную группу $AGL(d,q)$ как подгруппу группы $GL(d+1,q)$, сохраняющую $V$. 
Согласно пункту (а) выше, теорема~$\ref{thm:GL}'$ в размерности $d+1$ влечет 
те же оценки среднего размера стабилизатора применительно к действию $AGL(d,q)$ на $\Omega_{m}$ и $\Omega_{(m)}$. Однако 
по этим данным нельзя судить о среднем размере стабилизаторов действия $AGL(d,q)$, суженного на $V_m$ или $V_{(m)}$:
$m$-точечные конфигурации в гиперплоскости $V$ составляют малую часть всей совокупности сответствующих конфигураций в $\Omega$.



\subsection{Перечисление графов и гиперграфов}
\label{ssec:graphs}

Зафиксируем натуральное число $\ell\geq 2$. Пусть $X$ --- конечное множество, $|X|=n$;
будем считать $n$ большим параметром.
Положим $\Omega=X_\ell$ (множество $\ell$-элементных подможеств в $X$). 
Элементы множества $\Set{m}$ называются {\em помеченными однородными $\ell$-гиперграфами на $X$ с $m$ гиперребрами} \cite{Duchet1995}. При $\ell=2$ это помеченные графы с $m$ ребрами. 

На $\Omega$ естественным образом действует группа $Sym(X)$. Ее гомоморфный образ в полной группе подстановок $\Sym(\Omega)$ будем обозначать $S_n^\ell$. Как абстрактные группы, $\Sym(X)$ и $S_n^\ell$ (при $\ell\leq n-1$) изоморфны полной симметрической группе на $n$ индексах, и мы будем говорить о действии элемента (подстановки) $\sigma$ на $X$  или $\Omega$, имея в виду естественный изоморфизм $\Sym(X)\to S_n^\ell\subset \Sym(\Omega)$. Что касается естественного действия $S_n^\ell$ на  $\Omega_m$, то специальное обозначение для соответствующей группы подстановок --- гомоморфного образа  $S_n^\ell$ в $Sym(\Omega_m)$ --- нам не потребуется.
Однако говорить о цикловой структуре элемента $\sigma$ можно лишь зная, о действии на каком множестве идет речь. Особенно внимательно за этим надо следить в доказательстве теоремы~\ref{thm:graphs}\ ниже.

Рассмотрим задачу о перечислении {\em непомеченных}\ однородных гиперграфов
с $n$ вершинами и $m$ гиперребрами, т.е.\ $S_n^\ell$-орбит на $\Set{m}$.

Результаты ниже основаны на асимптотической свободе действия, но случаи $\ell\geq 3$ и $\ell=2$ требуют разных подходов: в первом, 
если не стремиться к оптимизации результата,
 достаточно оценить минимальную степень $G=S_n^\ell$ (действующей на $\Omega$) и воспользоваться теоремой~\ref{thm:number_of_orbits}, в то время как во втором мы по необходимости прибегнем к теореме~\ref{thm:number_of_orbits_2tier}. Все дело в соотношении величин минимальной степени
и логарифма порядка группы. При всех $\ell\geq 2$ минимальная степень имеет порядок $n^{\ell-1}$,
а $\log|G|\sim n\log n$, 
так что при $\ell\geq 3$ имеем $\mu\gg \log|G|$, а при $\ell=2$ соотношение обратное.

\begin{theorem}
\label{thm:hypergraphs}
Пусть фиксировано $\ell\geq 3$,
$$
|\Omega|={n\choose \ell},\qquad
|\Set{m}|={|\Omega|\choose m},\qquad
\overline{m}=\min(m, |\Omega|-m).
$$
Если $\overline{m},n\to\infty$ таким образом, что $\underline\lim\;\overline{m}(n^2\ln n)^{-1}>(2\ell)^{-1}$,
то
число $|\Set{m}/S_n^\ell|$ непомеченных однородных $\ell$-гиперграфов с $n$ вершинами и $m$ гиперребрами
имеет асимптотику 
\beq{est_mgraph}
 |\Set{m}/S_n^\ell|=\frac{|\Set{m}|}{n!}(1+o(1)).
\eeq
\end{theorem}

\proof
Достаточно рассмотреть случай $m\leq n$, т.е.\ $\overline m=m$. При достаточно больших $n$ имеем $m\geq \kappa n^2\ln n$, 
где $2\ell \kappa>1$. Соответственно $\lambda=m/|\Omega|\geq \kappa\ell! \, n^{2-\ell}\ln n$, и из (\ref{xi12})
следует, что $\ln\xi_1(\lambda)\leq -(\kappa\ell!-o(1)) \, n^{2-\ell}\ln n$.
 
Асимптотика (\ref{est_mgraph}) следует теперь из теоремы~\ref{thm:number_of_orbits} и оценки минимальной степени группы $S_n^{\ell}$
\beq{mindegSnell}
 \mu \geq  \frac{2n^{l-1}}{(l-1)!}(1-o(1)).
\eeq 
В самом деле, поскольку $\ln|S_n^\ell|=\ln n!\sim n\,\ln n$, логарифм правой части неравенства (\ref{general_stab_size-sets}) не превосходит
$n\ln n(1-2\ell\kappa+o(1))\to -\infty$.

\smallskip
Докажем оценку (\ref{mindegSnell}). Пусть $\sigma\in Sym(X)$, $|X|=n$, и $j_k=j_k(\sigma)$ --- число циклов длины $k$,
ср.\ \S~\ref{sec:est1}, п.~2.
Нам нужно оценить сверху число неподвижных точек действия $\sigma$ на $X_\ell$.
(Выведенная в \S~\ref{sec:est1} оценка (\ref{ub_alpham}) в данной ситуации бесполезна.) 
Искомое число $\alpha_m(\sigma)$ --- коэффициент при $t^\ell$ в произведении
$\prod_{k=1}^\ell p_k(t)$,
где $p_k(t)$ --- часть биномиального разложения $(1+t^k)^{j_k}$, включающая только члены степени
$\leq \ell$. Запишем $p_k(t)=\sum_{i=0}^\ell c_{k,i} t^i$. 
Ясно, что $c_{k,i}=0$, если $k$ не делит $i$, а
$$
 c_{k, kq}={j_k\choose q}\leq \frac{n^q}{q!}.  
$$
Имеем
$$
 \alpha_\ell(\sigma)=\sum 
 c_{1,i_1} c_{2,i_2}\dots c_{\ell, i_\ell},  
$$
суммирование ведется по таким индексным наборам $(i_1,i_2,\dots,i_\ell)$,
для которых $i_k=kq_k$ и $i_1+i_2+\dots+i_\ell = \ell$, считая ${j\choose q}=0$ при $q>j$. 
Число индексных наборов зависит только от $\ell$, но не от $n$.
Для произведения коэффициентов имеем оценку
$$
 c_{1,i_1} c_{2,i_2}\dots c_{\ell,i_\ell}\leq \frac{n^{Q}}{q_1!\dots q_\ell!},
$$
где   
$Q=\sum_{i=1}^\ell q_i
=\ell-q_2-2q_3-\dots-(\ell-1)q_\ell$.
Возможны случаи:

1) $Q=\ell$; так будет, если $q_1=\ell$, $q_2=q_3=\dots=0$.

2) $Q=\ell-1$; так будет, если $q_1=\ell-2$, $q_2=1$, $q_3=\dots=0$.

3) $Q\leq\ell-2$ для прочих индексных наборов.

\smallskip\noindent
Следовательно,
$$
 \alpha_\ell(\sigma) = c_{1,\ell} +c_{1,\ell-2} c_{2,2} +O(n^{\ell-2})
 ={j_1\choose \ell}+{j_1\choose \ell-2} \,{j_2\choose 1}\,+O(n^{\ell-2}).
$$
При $\sigma\neq e$ имеем $j_1\leq n-2$. Легко видеть, что при этом ограничении и условии $j_1+2j_2\leq n$ максимум суммы первых двух слагаемых достигается при $j_1=n-2$, $j_2=1$, и в этом случае второе слагаемое есть $O(n^{\ell-2})$.
Получаем нижнюю оценку степени подвижности подстановки $\sigma\in S_n^\ell$
$$
\begin{array}{rcl}
\beta_\ell(\sigma)&\geq&
\displaystyle
|\Omega|-{n-2\choose \ell} -O(n^{\ell-2})
= 
\\[2ex]
&=& \displaystyle
{n\choose \ell}-{n-2\choose \ell}-O(n^{\ell-2})
= 
\\[2ex]
&=& \displaystyle
\frac{2n^{\ell-1}}{(\ell-1)!}-O(n^{\ell-2}),
\end{array}
$$
что и требовалось.
\eop

\bigskip
Классический случай $\ell=2$ (перечисление {\em графов}) труднее в указанном ранее смысле.
Оценки, основанной только на минимальной степени, здесь недостаточно, но можно использовать технический прием <<изоляции плохого малого множества>>. 
Используя теорему~\ref{thm:number_of_orbits_2tier} с большим числом слагаемых, получим результат, близкий к неулучшаемому. 

\begin{theorem}
\label{thm:graphs}
Положим $\overline{m}=\min(m,|\Omega|-m)$.
Если $m$ и $n$  стремятся к бесконечности таким образом, что $\kappa:=\underline\lim\;\overline{m} (n\log n)^{-1}>1$,
то число $g_{n,m}=|\Set{m}/S_n^2|$ непомеченных графов с $n$ вершинами и $m$ ребрами
имеет асимптотику 
\beq{est_graph}
 g_{n,m}=\frac{|\Set{m}|}{n!}
 (1+o(1)).
\eeq
Здесь
$$
|\Omega|={n\choose 2},\qquad
|\Set{m}|={|\Omega|\choose m}.
$$
\end{theorem}

\noindent{\bf Замечание.}
Необходимое и достаточное условие асимптотической свободы действия $S_n^2$ установлено в работах \cite{Korshunov1970}, \cite{Wright1971}:
$\frac{2\overline{m}}{n}-\log n\to +\infty$ (в частности, {\it необходимо}\ условие $\kappa\geq 1/2$). В \cite{Korshunov1970} и \cite{Wright1972} найдена предельная форма распределения $g_{n,m}$ в критической области $m\sim \frac12 n\log n$.

\bigskip
\proof
Применим теорему~\ref{thm:number_of_orbits_2tier}, построив подходящую цепочку вложенных множеств, каждое из которых имеет вид
$B(r)=\{\sigma\in\Sym(X)\,|\,\alpha(\sigma)\geq n-r\}$. Прежде всего надо решить две задачи: 

А. Оценить сверху $|B(r)|$;

Б. Рассматривая $B(r)$ уже как подмножество группы $G_n^2$, оценить снизу минимальную степень его дополнения, $\mu(r)=\mu^{S_n^2\setminus B(r)}$. 

\medskip\noindent
{\it Решение задачи А. }
Пусть $\Sigma_r=\{(Y,\sigma)\,|\, Y\in X_r,\; \mathrm{supp}(\sigma)\subset Y\}$.
%
%
Отображение $\Sigma_r\to B(r)$, заданное формулой $(Y,\sigma)\mapsto\sigma$, сюръективно.
(Подстановка $\sigma\in B(r)$ имеет носитель размера $\leq r$ и есть образ пары $(Y,\sigma)$,
где $Y$ --- любое $r$-элементное надмножество носителя $\sigma$.)
При каждом $Y_0\in X_r$
слой $\{(Y,\sigma)\in\Sigma_r\,|\,Y=Y_0\}$ отождествляется с группой $Sym(Y_0)$, т.е.\ имеет размер $r!$. 
Следовательно,    
\beq{graph_probA}
 |B(r)|\leq {n\choose r}\,r! \leq n^r.
\eeq
 
\medskip\noindent
{\it Решение задачи Б. }
Рассмотрим действие подстановки $\sigma\in \Sym(X)\setminus B(r)$ параллельно на $X$ и на 
$\Omega=X_2$. Величины $j_1(\sigma)=\alpha(\sigma)$, $j_2(\sigma)$, $\beta(\sigma)$ относятся к действию на $X$, а величины $\alpha_2(\sigma)$ (ср.\ (\ref{alpham})) 
и $\beta_2(\sigma)=|X_2|-\alpha_2(\sigma)$ относятся к действию на $X_2$.

Неподвижные точки действия $\sigma$ на $\Omega$ бывают двух типов:

\smallskip
1) ${\alpha(\sigma)\choose 2}$ пар $(x,y)$, где $x$ и $y$ ---  неподвижные точки $\sigma$ на $X$;

\smallskip
2) $j_2(\sigma)$ пар $(x,y)$, соответствующих $2$-циклам $\sigma$ на $X$.

\smallskip\noindent
Итак,
$$
\alpha_2(\sigma)={j_1(\sigma)\choose 2}+j_2(\sigma).
$$
Всегда $j_1(\sigma)+2j_2(\sigma)\leq n$. Следовательно,  $j_2(\sigma)\leq \beta(\sigma)/2$ и 
$$
\beta_2(\sigma)=
{n\choose 2}-{j_1(\sigma)\choose 2}-j_2(\sigma)
\geq 
\frac{\beta(\sigma)(2n-\beta(\sigma)-2) }{2}.
$$
Если $\sigma\in \Sym(X)\setminus B(r)$, то $\beta(\sigma)>r$.  
Функция $r\mapsto r(2n-r-2)$ не убывает при $r\leq n-1$. Значит,
$$
\beta_2(\sigma) > \frac{r(2n-r-2)}{2}.
$$
Задача Б решена: получена оценка 
\beq{graph_probB}
\mu(r)=\mu^{S_n^2\setminus B(r)} > \frac{r(2n-r-2)}{2}.
\eeq

\medskip

Ниже предполагаем, что $m/n^2\to 0$ при $n\to \infty$; если это не так, то оценки проще. При сделанном предположении
имеем $\lambda=m/{n\choose 2}=o(1)$, и из (\ref{xi12}) следует, что $\xi_1(\lambda)=-2m/n^2(1+o(1))$.

Положим в теореме~\ref{thm:number_of_orbits_2tier} $B_i=B(r_i)$, $i=1,\dots, p-1$ (значения $p$ и $r_i$ подлежат определению)
и оценим правую часть (\ref{general_stab_size-sets3})  сверху суммой
$$
 3\sqrt{m}\left(|B(r_1)|\xi_1^\mu+\sum_{i=2}^{p-1}|B(r_i)|\xi_1^{\mu(r_{i-1})}+|G|\xi_1^{\mu(r_{p-1})}\right)=:\sum_{i=1}^{p} K_p.
$$
Напомним: символ $n$ в (\ref{general_stab_size-sets3}) теперь соответствует $|\Omega|=\frac12 n^2 (1+o(1))$.

Минимальная степень группы $S_n^{2}$ равна $\mu=2n-4$.
(Ср.\ доказательство теоремы~\ref{thm:hypergraphs} или \cite{DixonMortimer1996}, Exercise 5.3.4.) 
С учетом этого факта и неравенств (\ref{graph_probA}), (\ref{graph_probB}) имеем 
$$
\begin{array}{l}
\dst
 \ln K_1 \leq\ln\left(3{m}^{1/2}|B_1|\xi_1^\mu\right)
 \leq \left(\frac{1}{2}+r_1\right)\ln n-\frac{4m}{n}(1-o(1)),
\\[3ex]\dst
\ln \sum_{i=2}^{p-1} K_i\leq \ln p+\frac{\ln n}{2}+\max_{2\leq i\leq p-1}\left(r_i\ln n- \frac{r_{i-1}(2n-r_{i-1})m}{n^2}(1-o(1))
\right),
 \\[3ex]
 \dst
 \ln K_p \leq\ln\left(3{m}^{1/2}|G|\xi_1^{\mu(r_{p-1})}\right)
 \leq n\ln n-\frac{r_{p-1}(2n-r_{p-1})m}{n^2}(1-o(1)).
\end{array}
$$
Покажем, что при условии $m/n\geq (\kappa-o(1))\ln n$, $\kappa>1$, можно выбрать $r_i$ ($i=1,\dots,p-1$) так, чтобы 
правые часть этих неравенств стремились к $-\infty$. 
В первом неравенстве достаточно положить $r_1=2$. 
Далее последовательно выбираем $2<r_2<r_3<\dots<r_p<n$  согласно правилу: $r_i$ есть наибольшее целое, удовлетворяющее неравенству
$r_{i}+1/2<\kappa r_{i-1} (2 -r_{i-1}/n)$. Поскольку $\kappa>1$, то при достаточно большом $n$ имеем $r_2\geq 3$, и последовательность
$\{r_i\}$ монотонно возрастает. 
Легко видеть, что при $n\geq n_0(\kappa)$ построенная последовательность растет не медленнее геометрической прогрессии, 
$r_i\geq Ch^i$, где $C>0$ и $h>1$ зависят от $\kappa$, но не от $n$. Следовательно, $p=O(\ln n)$. По выбору $\{r_i\}$, правая часть второго неравенства стремится к $-\infty$. 
Последний член последовательности, $r_{p-1}$, по построению удовлетворяет неравенству 
$n+1/2\leq \kappa r_{p-1} (2 -r_{p-1}/n)$. Поэтому 
$$
\begin{array}{c}
\dst
 n\ln n-\frac{r_{p-1}(2n-r_{p-1})m}{n^2}(1-o(1))
\leq n\ln n -\left(n+\frac{1}{2}-o(1)\right)\ln n,
 \end{array}
$$
и правая часть третьего неравенства также стремится к $-\infty$. Теорема доказана.
\eop

\bigskip\noindent
{\bf Замечание.}
Теорема~\ref{thm:number_of_orbits_2tier} в двучленной версии (\ref{general_stab_size-sets2}) позволяет
проще доказать асимптотику (\ref{est_graph}) в предположении $\underline\lim\;m (n^{3/2}\log n)^{-1}>2^{-3/2}$.


\subsection{Оценки для общих групп подстановок}
\label{ssec:gengroups}

В теории групп подстановок известны универсальные оценки минимальной степени $\mu$ снизу 
и порядка группы $|G|$ сверху в терминах степени $n=|\Omega|$ (размера пространства действия).
Это как раз входная информация для теоремы~\ref{thm:number_of_orbits}. Для простоты в этом разделе предполагаем, 
что $\lambda=m/n$ и $\eta=m/(n+m)$ отделены от $0$ и $1$.
В пункте 2 (с.~\pageref{thm:Babai}) покажем, что теорема~\ref{thm:number_of_orbits} дает  
содержательное усиление тривиальной верхней оценки в (\ref{triv_num_orbits})
для действий больших групп на множествах больших (мульти)множеств в случае, когда базовая группа подстановок примитивна или 2-транзитивна, а в остальном --- почти произвольна. Исключение составляют полная симметрическая и знакопеременная группы. 
В пункте~1 для большей полноты картины анализируем качество оценок теорем~\ref{thm:number_of_orbits} 
и \ref{thm:number_of_orbits_2tier} в случае действия симметрической группы (для знакопеременной группы анализ полностью аналогичен) на подмножествах. 

\medskip\noindent{\bf 1. }
Действие полной симметрической группы $G=Sym(\Omega)$.
на $\Set{m}$ 
наиболее далеко от свободного. Разумеется, интерес здесь представляет не заключение теоремы~\ref{thm:number_of_orbits}
, а лишь сравнение даваемых ею оценок с истинным числом орбит по порядку величины: это позволит судить о том, насколько грубы предсказания теоремы в максимально неблагоприятной ситуации.

\smallskip
Действие $Sym(\Omega)$ на $\Set{m}$ транзитивно при всех $m\leq n$, следовательно, $|\Set{m}/G|=1$ и 
$$
\frac{1}{|G|}\avStab{G}{m}=\frac{1}{|\Set{m}|}\sim e^{-nh(\lambda)} n^{1/2}\cdot \mathrm{const}(\lambda)
$$
(ср.\ (\ref{stirling-binomial-1})).
 
Минимальная степень полной симметрической группы, $\mu=2$, не растет с ростом порядка группы. Неравенство  (\ref{general_stab_size-sets}) дает лишь (неинформативную) оценку  
$$
\frac{1}{|G|}\avStab{G}{m}\leq n^{1/2}\cdot\mathrm{const}(\lambda).
$$
Мы видим, что оценка (\ref{general_stab_size-sets}) завышена, по порядку величины, в
$$
 e^{n h(\lambda)}\sim |G|^{h(\lambda)/\ln\ln|G|} 
$$
раз. (По формуле Стирлинга $n\sim \ln|G|/\ln\ln|G|$.)

С помощью теоремы~\ref{thm:number_of_orbits_2tier} можно сократить указанный разрыв. Для простоты воспользуемся ее двучленным вариантом (\ref{general_stab_size-sets2}).
Выделим <<плохое>> множество $B=B(r)\subset G$ как в доказательстве теоремы~\ref{thm:graphs}; пусть $r=tn$, где $t<1$ ---
параметр. Имеем $\mu^{G\setminus B}>r$ и, согласно (\ref{graph_probA}), $|B|\leq n^{tn}$.  В правой части неравенства (\ref{general_stab_size-sets2}) первое слагаемое асимптотически доминирует, однако наличие множителя $\xi_1^{\mu^*}$, где $\mu^*=tn$, в отличие от $\mu$, растет вместе с $n$, приводит к оценке вида   
\beq{symn_improved}
\frac{1}{|G|}\avStab{G}{m}\leq e^{tn\ln\xi_1(\lambda)} n^{1/2}\cdot\mathrm{const}(\lambda),
\eeq
которая, заметим, уже информативна, т.е.\ лучше тривиальной.
Эта оценка завышена, по порядку величины, в
$$
 e^{n q(t,\lambda)}\sim |G|^{q(t,\lambda)/\ln\ln|G|} 
$$
раз, где
$$
q(t,\lambda)=h(\lambda)+t\ln\xi_1(\lambda).
$$
<<Коэффициент уменьшения ошибки>> в логарифмической шкале за счет использования неравенства (\ref{general_stab_size-sets2}) 
по сравнению с (\ref{general_stab_size-sets}) равен
$$
 \frac{q(t,\lambda)}{ h(\lambda)}=1+ t\frac{\ln \xi_1(\lambda)}{ h(\lambda)}.
$$
При $\lambda=1/2$ (где достигается 
$\min \ln\xi_1(\lambda)/h(\lambda)={-1/2}$), 
эту величину можно сделать сколь угодно близкой к $1/2$, беря $t$ сколь угодно близким к 1.
(Условно говоря, еще одно улучшение аналогичной силы --- и получилась бы оценка иного качества.)  
Добиться существенно большего теорема~\ref{thm:number_of_orbits_2tier} не позволяет даже при использовании длинной цепочки вложенных множеств.
В самом деле, рассмотрим слагаемое с индексом $p$ в правой части (\ref{general_stab_size-sets3}). Наибольшее возможное значение показателя
$\mu_{p-1}$ равно $n$, а наименьшее возможное значение $|G|-|B_{p-2}|$ равно числу подстановок, не имеющих неподвижных точек, которое, как
известно, имеет асимптотику $e^{-1}|G|$ при $n\to\infty$. Отсюда следует, что 
невозможно получить оценку, асимптотически лучшую, чем (\ref{symn_improved})  с $t=1$.

\bigskip\noindent
{\bf 2. } В случае транзитивных групп, отличных от полной симметрической и знакопеременной,
теорема~\ref{thm:number_of_orbits} ведет к содержательным следствиям.  

Группа подстановок $G$ на множестве $\Omega$ {\it примитивна}, 
если ни при каком $k\geq 2$ не существует орбиты действия $G$ на $\Set{k}$, состоящей из попарно непересекающихся подмножеств. (Обратно, если такая орбита существует, то ее члены называются {\it блоками импримитивности}.)\@ Очевидно, что примитивная группа транзитивна. Легко видеть также, что {\it $2$-транзитивное}\ действие $G$ (такое, что индуцированное действие на $\Set{2}$ транзитивно) всегда примитивно  \cite[Theorem 9.6]
{Wielandt}, но обратная импликация неверна. 
Примитивная, но не $2$-транзитивная группа называется {\it унипримитивной}. Таким образом, любая транзитивная группа либо $2$-транзитивна, либо унипримитивна.

В этом пункте опираемся на результаты Бабаи (усиление теоремы Жордана 1871 г.) \cite{Babai1981}, \cite[Theorem~5.3A]{DixonMortimer1996}, \cite{Zemlyachenko1985}, \cite[p.~2054]{LPWZ1995}, Бохерта \cite{Bochert1897}, \cite{HerzogPraeger}, \cite[\S~15]{Wielandt} и Пибера 
\cite{Pyber1993}, \cite[Theorem~5.6A]{DixonMortimer1996}. 
Как и раньше, $\mu$ обозначает минимальную степень $G$.

\medskip\noindent
{\bf Теорема Бабаи. }
\label{thm:Babai} 
{\it Для унипримитивной группы $G$ справедливы неравенства
\beq{Babai_ineq_mindeg}
 \mu > \sqrt{n}/2
\eeq
и
\beq{Babai_ineq_order}
|G|<\exp\big(4\sqrt{n}(\ln n)^2\big).
\eeq
}

\medskip\noindent
{\bf Комбинация теорем Бохерта и Пибера. } 
\label{thm:Bochert} 
{\it Для $2$-транзитивной группы $G\not\supset Alt(\Omega)$ имеет место  неравенство Бохерта
\beq{Bochert_ineq_mindeg}
 \mu > \frac{n}3-\frac{2}{3}\sqrt{n} 
\eeq
и неравенство Пибера
\beq{Pyber_ineq_order}
|G|<\exp\big(72(\ln n)^3\big).
\eeq
}

%

%
Подставляя неравенства Бабаи, Бохерта и Пибера в теорему~\ref{thm:number_of_orbits}, получим усиление тривиальной верхней оценки среднего размера стабилизатора индуцированного действия $G$ на множествах $\Set{m}$ и $\MSet{m}$. Усиление охватывает все достаточно большие транзитивные группы, за исключением $Sym$ и $Alt$.

В формулировке теоремы пишем $\xi=\xi_1(\lambda)$ и $\avStab{}{}=\avStab{G}{m}$ для действия на подмножествах, соответственно,
$\xi=\xi_2(\eta)$ и $\avStab{}{}=\avStab{G}{(m)}$ для действия на мультимножествах; функции $\xi_1(\lambda)$ и $\xi_2(\eta)$ определены в {\rm (\ref{xi12})}.

\begin{theorem}
\label{numorbits_primitive_groups}
Для любых $\delta\in (0,1/8)$ и $\eps\in (0,\,1/2)$ существует $n_0=n_0(\delta,\eps)$ такое, что для транзитивных групп подстановок на множестве 
 $\Omega$, $|\Omega|=n>n_0$ верны следующие утверждения.
 
\smallskip {\rm (а)}
Если группа $G$ унипримитивна, то при естественном действии $G$ на подмножествах множества $\Omega$ мощности $m\in (n\eps,\, n(1-\eps))$ и на мультимножествах веса $m\in (n\eps/(1-\eps) ,\, n(1-\eps)/\eps)$ существует хотя бы одна орбита $\mathcal{O}$ размера 
\beq{est_orbitsize_uniprimitive}
|\mathcal{O}|\geq |G|^{A/(\ln\ln |G|)^2},
\qquad
A=\left(\frac{1}{8}-\delta\right)|\ln\xi|.
\eeq

\smallskip {\rm (б)}
Если группа $G$ $2$-транзитивна и отлична от $Sym(\Omega)$ и $Alt(\Omega)$, то 
доля регулярных орбит среди всех орбит естественного действия группы $G$ на подмножествах $\Omega$ мощности $m\in (n\eps,\, n(1-\eps))$ и на мультимножествах веса $m\in (n\eps/(1-\eps) ,\, n(1-\eps)/\eps)$
есть $1-o(1)$ при $n\to\infty$. Более точно, отклонение от $1$ среднего размера стабилизатора в обоих случаях 
 допускает оценку сверху:
\beq{est_2transitive}
\avStab{}{}-1
\leq
\xi^{n/3-(2/3+\delta)\sqrt{n}}.
\eeq
\end{theorem}

\proof
(а) 
Неравенство (\ref{Babai_ineq_mindeg}) и теорема~\ref{thm:number_of_orbits} дают
$$
\frac{\avStab{}{}}{|G|} 
  \leq\mathrm{const}(\eps)\cdot\sqrt{n}\cdot\xi^{\sqrt{n}/2}=\exp\left(\frac{\ln\xi}{2}\sqrt{n}(1+o(1))\right).
$$
Пусть $S$ обозначает $\Set{m}$ или $\MSet{m}$. 
Очевидно неравенство
$$
 |S|\leq |S/G|\cdot \max_{\mathcal{O}\in S/G} |\mathcal{O}|.
$$
Следовательно, существует хотя бы одна орбита $\mathcal{O}$, для которой
$$
|\mathcal{O}|\geq \frac{|S|}{|S/G|}=\frac{|G|}{\avStab{}{}}
\geq \exp\left(\frac{|\ln\xi|}{2}\sqrt{n}(1-o(1))\right).
$$
Неравенство (\ref{Babai_ineq_order}) влечет, как легко видеть,
$$
 \sqrt{n}\geq \frac{\ln|G|}{4(\ln\ln|G|)^2}(1-o(1)).
$$
Комбинируя два последних неравенства, получаем (\ref{est_orbitsize_uniprimitive}).

\medskip
(б) Для $2$-транзитивной группы $G\not\supset Alt(\Omega)$ по теореме~\ref{thm:number_of_orbits} с учетом (\ref{Bochert_ineq_mindeg}) имеем
$$
\frac{\avStab{G}{m}-1}{|G|}
  \leq \xi^{n/3-2/3\sqrt{n}+O(\ln n)},
$$
а согласно (\ref{Pyber_ineq_order}), $|G|=\exp(o(\sqrt{n}))$. 
Так получается оценка (\ref{est_2transitive}).
\eop


\bigskip\noindent
{\bf Замечания. }
1. Использованные теоретико-групповые результаты не зависят от классификации простых конечных групп.
Существуют более точные оценки для широких классов групп с описанными исключениями --- см., например, \cite{LiebeckShalev01,Maroti02}. 

\smallskip
2. Действие полной аффинной группы, рассмотренное в разд.~\ref{ssec:affine}, $2$-транзитивно, поэтому версия теоремы~\ref{thm:GL} с худшими остаточными членами является следствием теоремы~\ref{numorbits_primitive_groups}(б). Напротив, действие группы $G=Sym(V)$ перенумерации вершин графа на $\Omega=V_2$,
унипримитивно, и формулировка теоремы~\ref{thm:graphs} намного сильнее, чем общей теоремы~\ref{numorbits_primitive_groups}(а). Это неудивительно: доказательство теоремы~\ref{thm:graphs} опирается не только на оценку минимальной степени группы, но и на оценки минимальных степеней
ее подмножеств $B_i$   

\smallskip
3. Интересно было бы проанализировать применимость неравенства (\ref{est_2transitive}) с целью воспроизведения (не-асимптотического) 
результата работы \cite{SiemonsZalesski02} о существовании регулярных орбит действия циклических подгрупп известных простых групп,
действующих 2-транзитивно.

\smallskip
4. Было бы полезно иметь обобщения оценок типа (\ref{Babai_ineq_mindeg})--(\ref{Pyber_ineq_order}) в направлении,
подсказываемом теоремой~\ref{thm:number_of_orbits_2tier}. Помимо ранее цитированных работ, 
содержащих результаты и ссылки, которые могут при этом оказаться полезными, укажем на обзор на русском языке \cite{GlukhovZubov00}.  

\bigskip
Автор признателен М.В.~Кочетову за предложение принять участие в проекте \cite{KPS13} и обсуждения, в результате которых оформилась
тема настоящей работы. Также автор благодарен рецензенту за содержательные замечания, способствовавшие усилению первоначальных формулировок. 


\pagebreak


\addcontentsline{toc}{section}{Список литературы}
\begin{thebibliography}{30}





%


%

\bibitem{Harary-Palmer} 
Харари Ф., Палмер Э. {\it Перечисление графов. } Мир, Москва, 1977.

\bibitem{FordUhlenbeck1957-4} Ford G. W., Uhlenbeck G. E.
Combinatorial problems in the theory of graphs IV, {\it Proc. Nat. Acad. Sci. USA}\ (1957) {\bf 43},
163--167.

\bibitem{Hoffman-Welch1968}
Hoffman F., Welch L. Totally variant sets in finite groups and vector spaces.
{\it Canadian J. Math.}\ (1968) {\bf 20}, 701--710.


\bibitem{Strazdins1997}
Strazdins I. Universal affine classification of boolean functions. {\it Acta Applic. Math.}\ (1997) {\bf 46}, 147--167.

\bibitem{Wild00} Wild M. The asymptotic number of inequivalent binary codes and nonisomorphic binary matroids. {\it Finite Fields and Their Appl.}\ (2000) {\bf 6}, 192--202.

\bibitem{Wild05} Wild M. The asymptotic number of binary codes and binary matroids.
{\it SIAM J. Discrete Math.}\ (2005) {\bf 19}, \No~3, 691--699.

\bibitem{Hou07a}
Xiang-Dong Hou, On the asymptotic number of non-equivalent binary linear codes, 
{\it Finite Fields and Their Appl.}\ (2007) {\bf 13}, 318--326.

\bibitem{Hou07b}
Xiang-Dong Hou, On the asymptotic number of inequivalent binary self-dual codes. 
{\it J.~Combinat.~Theory Ser.~A}\ (2007) {\bf 114}, 522--544. 

\bibitem{Hou09} 
Xiang-Dong Hou, Asymptotic numbers of non-equivalent codes in three notions
of equivalence. {\it Linear and Multilinear Algebra}\ (2009) {\bf 57}, \No~2, 111--122. 

\bibitem{Lax04} Lax R. F.
On the character of $S_n$ acting on subspaces of $\FF_n$.
{\it Finite Fields and Their Appl.}\ (2004) {\bf 10}, 315--322.

%


\bibitem{Cameron1977} Cameron  P. J. Permutation groups on unordered sets. 
In: {\it Higher Combinatorics}\ (M.~Aigner, Ed.). Proc. NATO Adv. Study Inst. Series {\bf 31},  
D.~Reidel Publ. Co., Amsterdam, 1978, pp.~217--239. 
%
Рус.\ перевод: Камерон П.~Дж. Группы подстановок на неупорядоченных множествах. В кн.: {\it Проблемы комбинаторного анализа}\ (ред. Рыбников К. А.). Мир, Москва, 1980, с.~151--178.

\bibitem{Cameron1999} Cameron P.~J. {\it Permutation groups. } Cambridge Univ. Press, Cambridge, 1999.

\bibitem{Siemons1984} Siemons J. Permutation groups on unordered sets I. {\it Arch. Math.}\ (1984) {\bf 43}, 483--487.

\bibitem{LW1965}
Livingstone D., Wagner  A.
Transitivity of finite permutation groups on unordered sets.
{\it Math.\ Zeitschr.}\ (1965) {\bf 90}, 393--403.

\bibitem{Laborde1985} 
Laborde J.-M. Regularisation numerique d'orbites, {\it Discrete Math.} (1985)
{\bf 53}, 151--155. 

\bibitem{KPS13} Kochetov M., Parsons N., Sadov S., Counting fine gradings on matrix algebras and on classical simple Lie algebras.
{\it Int. J. of Algebra and Computations}\ (2013) 
{\bf 23}, \No~7, 1755--1781.

\bibitem{Korshunov1970} Коршунов, А.Д. О мощности некоторых классов графов.
{\it Доклады АН СССР} (1970) {\bf 193}, \No~6, 1230--1233.

\bibitem{Wright1971} Wright E. M.
Graphs on unlabelled nodes with a given number of edges. {\it Acta Math. } (1971) {\bf 168},
1--9.

\bibitem{Wright1972} Wright E. M.
The number of unlabelled graphs with many nodes and edges. {\it Bull. AMS} (1972) {\bf 78:6},
1032--1034.

\bibitem{Cameron1986} Cameron P.~J. Regular orbits of permutation groups on the
power set, {\it Discrete Math.} (1986), {\bf 62}, 307--309. 

\bibitem{AbramowitzStegun} {\it Справочник по специальным функциям} (ред.\ Абрамовиц М., Стиган И.)
Наука, Москва, 1979.

\bibitem{Wielandt} Wielandt H. {\it Finite permutation groups. }
 Academic Press, New York-London, 1964.

\bibitem{DixonMortimer1996} Dixon J. D., Mortimer B. {\it Permutation groups. }
Springer-Verlag, New York, 1996.

\bibitem{Duchet1995} Duchet P., Hypergraphs. In: {\it Handbook of Combinatorics}\ (Graham R., Gr\"otschel M., Lov\'asz L., Eds.). Elsevier Science, 
Amsterdam-Lausanne-New York-Oxford-Shannon-Tokyo, 1995, pp.~381--432.


\bibitem{Babai1981} Babai L., On the Order of Uniprimitive Permutation Groups,
{\it Annals of Math.}\ (1981) {\bf 113}, \No~3, 553--568.


\bibitem{Zemlyachenko1985} Земляченко В.Н., Корнеенко Н.М., Тышкевич Р.И.,  Проблема изоморфизма графов. В кн.: {\it Теория сложности вычислений.~I. Зап. научн. сем. ЛОМИ}\ (1982) {\bf 118}, с.~83--158.

\bibitem{LPWZ1995} Lov\'asz L., Pyber L., Welsh D. J. A., Ziegler G. M.  Combinatorics in Pure Mathematics. In: {\it Handbook of Combinatorics}\ (Graham R., Gr\"otschel M., Lov\'asz L., Eds.). 
Elsevier Science, 
Amsterdam-Lausanne-New York-Oxford-Shannon-Tokyo, 1995, pp.~2039--2082.

\bibitem{Bochert1897} Bochert A. Ueber die Classe der Transitiven Substitutionengruppen II.
{\it Math. Ann.}\ (1897) {\bf 49}, 133--144.

\bibitem{HerzogPraeger} Herzog M., Praeger C. E. Minimal degree of primitive permutation groups.
In: {\it Combinatorial Mathematics IV, Adelaide 1975}\  (Casse L. R. A. and Wallis W. D., Eds.). Lecture Notes in Math. {\bf 560}, Springer-Verlag, Berlin-Heidelberg-New York, pp.~116--122.

\bibitem{Pyber1993} Pyber  L. On the Orders of Doubly Transitive Permutation Groups, Elementary Estimates. {\it J. Combinat. Theory, Ser. A}\ (1993), {\bf 62}, 361--366.

\bibitem{LiebeckShalev01}
Liebeck M. W., Shalev A., Bases of primitive permutation groups. In: {\it Groups, Combinatorics and Geometry: Durham, 2001.}\ (Ivanov A. A., Liebeck M. W., Saxl J., Eds.). World Scientific, New Jersey-London-Singapore-Hong Kong, 2003, pp.~147--154.

\bibitem{Maroti02}
Mar\'oti A.
On The Orders of Primitive Groups.  {\it J. of Algebra}\ (2002), {\bf 258}, \No~2, 631--640. 

\bibitem{SiemonsZalesski02}
Siemons J., Zalesski\u{\i} A.   
Regular orbits of cyclic subgroups in permutation representations of certain simple groups.
{\it J. of Algebra} (2002), {\bf 256}, \No~2, 611--625.

\bibitem{GlukhovZubov00} Глухов М. М., Зубов А. Ю., О длинах симметрических и знакопеременных групп подстановок в различных системах образующих (обзор). В кн.: {\it Мат. вопр. кибернетики}, {\bf 8}, 
Наука, Москва, 2000, с.~3--30.

\end{thebibliography}
\end{document}